\documentclass[a4paper,reqno,10pt,oneside]{amsart}
\usepackage{amsmath,geometry}
\usepackage{graphicx}
\usepackage{mathrsfs}
\usepackage{amssymb,amsmath, amsfonts, amsthm}
\usepackage{latexsym}
\usepackage{color} 
\usepackage[dvips]{epsfig}
\usepackage{graphicx}

\allowdisplaybreaks[1]

\numberwithin{equation}{section}

\renewcommand{\O}{\operatorname{O}}

\renewcommand{\(}{\left(}
\renewcommand{\)}{\right)}
\renewcommand{\[}{\left[}
\renewcommand{\]}{\right]}

\newlength{\defbaselineskip}
\newtheorem{theorem}{Theorem}[section]

\newtheorem{lemma}[theorem]{Lemma}
\newtheorem{remark}[theorem]{Remark}

\renewcommand{\le}{\leqslant}
\renewcommand{\ge}{\geqslant}

\renewcommand{\l }{\lambda}
\newcommand{\n }{\nabla }

\renewcommand{\O}{\Omega}

\renewcommand{\L}{\Lambda}

\newcommand{\Sl}{\Sigma_\lambda}

\newcommand{\beq}{\begin{equation}}
\newcommand{\eeq}{\end{equation}}
\newcommand{\beqs}{\begin{equation*}}
\newcommand{\eeqs}{\end{equation*}}
\newcommand{\beqn}{\begin{eqnarray}}
\newcommand{\eeqn}{\end{eqnarray}}
\newcommand{\beqns}{\begin{eqnarray*}}
\newcommand{\eeqns}{\end{eqnarray*}}
\newcommand{\bdoc}{\begin{document}}
\newcommand{\edoc}{\end{document}}
\newcommand{\be}{\begin{enumerate}}
\newcommand{\ee}{\end{enumerate}}
\newcommand{\bdescr}{\begin{description}}
\newcommand{\edescr}{\end{description}}
\newcommand{\ba}{\begin{array}}
\newcommand{\ea}{\end{array}}
\newcommand{\intR}{\int_{\mathbb R^N}}
\newcommand{\R}{\mathbb R^N}
\newcommand{\e}{\varepsilon}
 \renewcommand{\(}{\left(}
\renewcommand{\)}{\right)}
\renewcommand{\[}{\left[}
\renewcommand{\]}{\right]}
\newenvironment{Proof}{\noindent{\bf Proof}}{\hfill$\Box$\\[2mm]}

\def\red{\color{red}}
\def\bl{\color{blue}}
\def\bk{\color{black}}

\begin{document}

\title[Radial symmetry for a quasilinear elliptic equation]{Radial symmetry for a quasilinear elliptic equation with a critical Sobolev growth and Hardy potential}

\author{Francescantonio Oliva, Berardino Sciunzi and Giusi Vaira}
\address{Francescantonio Oliva\\ Dipartimento di Scienze di Base e Applicate per l'Ingegneria, Sapienza Universit\`a di Roma\\Via Scarpa 16, 00161 Roma, Italy}
\email{francesco.oliva@sbai.uniroma1.it}
\address{Berardino Sciunzi\\ Dipartimento di Matematica e Informatica, UNICAL,, Ponte Pietro Bucci 31B, 87036 Arcavacata di Rende, Cosenza, Italy}
\email{sciunzi@mat.unical.it}
\address{Giusi Vaira\\ Dipartimento di Matematica e Fisica, Universit\`a degli studi della Campania ``Luigi Vanvitelli"\\Viale Lincoln 5, 81100 Caserta, Italy}
\email{giusi.vaira@unicampania.it}
\maketitle
\begin{abstract}
We consider weak positive solutions to the critical $p$-Laplace equation with Hardy potential in $\mathbb R^N$ $$-\Delta_p u -\frac{\gamma}{|x|^p} u^{p-1}=u^{p^*-1}$$ where $1<p<N$, $0\le \gamma <\left(\frac{N-p}{p}\right)^p$ and $p^*=\frac{Np}{N-p}$.\\ The main result is to show that all the solutions in $\mathcal D^{1, p}(\mathbb R^N)$ are radial and radially decreasing about the origin.

\end{abstract}
\section{Introduction and statement of the main result}
We study the doubly critical problem
 \begin{equation}\label{pb}
\left\{\begin{aligned}
&-\Delta_p u -\frac{\gamma}{|x|^p} u^{p-1}=u^{p^*-1}\qquad&\mbox{in}\,\,\mathbb R^N\\
&u>0\qquad & \mbox{in}\,\, \mathbb R^N\\
&u\in \mathcal D^{1, p}(\mathbb R^N)
\end{aligned}
\right.
\end{equation}
where $\Delta_p u := {\rm div}(|\nabla u|^{p-2} \nabla u)$ is the $p$-Laplacian operator with $1<p<N$,  $0\le \gamma<\gamma_p:=\left(\frac{N-p}{p}\right)^p$ and $p^*:=\frac{Np}{N-p}$ is the critical exponent for the Sobolev embedding.
Here $\mathcal D^{1, p}(\mathbb R^N)$ denotes the completition of $C^\infty_0(\mathbb R^N)$, the space of smooth functions with compact support, with respect to the norm $$\|u\|:=\left(\intR|\nabla u|^p\right)^{\frac 1p}.$$
By standard regularity theory, see \cite{Di,T}, it follows that solutions to \eqref{pb} are of class $C^{1,\alpha}$ far from the origin.\\

\noindent We address the study of the classification of positive solutions to \eqref{pb}. As we shall discuss later on, this is a crucial issue since problem \eqref{pb} naturally appears in the study of  $p$-Hardy-Sobolev inequalities as well as it appears as a limiting problem in many applications.
Our main effort is to show that all the positive solutions to \eqref{pb} are radial (and radially decreasing) about the origin. Once the radial symmetry of the solution is proved it is easy to derive the associated ordinary differential equation fulfilled by the solution $u=u(r)$. The
classification result  reduces therefore to an ODE analysis that has been already carried out in \cite{AFP} where the radial symmetry of the solutions was an assumption.\\

\noindent Let us start discussing the simpler case  $\gamma=0$. In this case the problem reduces to the following critical one
\begin{equation}\label{crit}
\left\{\begin{aligned}
&-\Delta_p u =u^{p^*-1}\quad&\mbox{in}\,\mathbb R^N,\\
&u>0\quad & \mbox{in}\, \mathbb R^N,\\
&u\in \mathcal D^{1, p}(\mathbb R^N).
\end{aligned}
\right.
\end{equation}
For such a problem a huge literature is available and the classification of positive weak solutions of \eqref{crit} is well understood. Indeed, for $\delta>0$ and $x_0\in\mathbb R^N$, an explicit family of solutions to \eqref{crit} is given by
\begin{equation}\label{solcrit}
V_{\delta, x_0}(x):=\left(\frac{\delta^{\frac{1}{p-1}}\alpha_{N, p}}{\delta^{\frac{p}{p-1}}+|x-x_0|^{\frac{p}{p-1}}}\right)^{\frac{N-p}{p}},
\end{equation}
where $\alpha_{N, p}:=N^{\frac 1 p}\left(\frac{N-p}{p-1}\right)^{\frac{p-1}{p}}$.
The family of functions given by \eqref{solcrit} are the minimizers to
\begin{equation}
\label{sp}
\mathcal{S}_p:=\inf_{\substack{\varphi\in \mathcal D^{1, p}(\mathbb R^N) \\ \varphi\neq 0}}\displaystyle\frac{\displaystyle\intR |\n\varphi|^p}{\left(\displaystyle\intR\varphi^{p*}\right)^{\frac{p}{p*}}}
\end{equation}
and the classification of the minimizers (see \cite{Talenti}) follows via symmetrization arguments. Note that   such a technique can be applied in the same way both in the semilinear case $p=2$ and in the quasilinear case $1<p<\infty$. \\

Furthermore, if we restrict the attention to the class of \emph{radial} solutions, then the analysis carried out in
\cite{GV} shows that all the \emph{regular radial} solutions to \eqref{crit} are given by \eqref{solcrit}.\\

For $p=2$ all the solutions to the equation are classified by \eqref{solcrit} as a consequence of the results in \cite{CGS} where  the Kelvin transform is strongly  exploited. A Kelvin type transformation is not applicable for the quasilinear case and this fact causes that a different proof is needed.
When no a priori assumption are imposed, the classification of all the positive solutions to \eqref{crit} (showing that all the solutions to \eqref{crit} are given by \eqref{solcrit}) has been in fact an open and challenging  problem recently solved in \cite{DMMS,S,Vetois} (see also   \cite{DPR, DR}). The techniques used are mainly based on a fine asymptotic analysis at infinity and refined versions of the \emph{moving plane procedure}, see \cite{GNN,Ser}.\\

Let us now turn to the case  $0< \gamma<\gamma_p$ but in the case $p=2$ so that $\gamma_2$ is the best constant in the Hardy-Sobolev inequality for $p=2$. For
$$\mathcal{S}_{2, \gamma}=\inf_{\substack{\varphi\in \mathcal D^{1, 2}(\mathbb R^N) \\ \varphi\neq 0}}\frac{\displaystyle\intR \left(|\n \varphi|^2-\frac{\gamma}{|x|^2}\varphi^2\right)}{\left(\displaystyle\intR |\varphi|^{2^*}\right)^{\frac{2}{2^*}}},$$
 it is known that $\mathcal{S}_{2, \gamma}$ is attained and extremals for $S_{2, \gamma}$ have the form (up to a multiplicative constant)
\begin{equation}\label{bubblebis}
U_\delta(x)=\delta^{-\frac{N-2}{2}} U\(\frac{x}{\delta}\)=\frac{\alpha_N \delta^{\Gamma  }}{|x|^{\beta_-}(\delta^{\frac{4\Gamma}{N-2}}+|x|^{\frac{4\Gamma}{N-2}})^{\frac{N-2}{2}}},\quad  \delta>0, \end{equation}
where
\begin{equation*}\label{bubble0bis} U (x)=\frac{\alpha_N}{|x|^{\beta^-}(1+|x|^{\frac{4\Gamma}{N-2}})^{\frac{N-2}{2}}}=
\frac{\alpha_N}{\(|x|^{\frac2{N-2}\beta^-}+|x|^{\frac2{N-2}\beta^+}\)^{\frac{N-2}2}}
\end{equation*}
with
\begin{equation*}\label{betabis}
\Gamma=\sqrt{\frac{(N-2)^2}{4}-\gamma},\quad \beta_\pm=\frac{N-2}{2}\pm \Gamma, \quad \alpha_N=\[\frac{4 \Gamma^2 N}{N-2}\]^{\frac{N-2}{4}},
\end{equation*}
see \cite{CaWa,ChCh,Terracini}. Moreover \eqref{bubblebis} gives all the solutions of the problem \eqref{pb} for $p=2$ and $\gamma\in (0, \gamma_2)$
and this has been proved in the celebrated paper \cite{Terracini}. In the case $p=2$ it is also known that when $\gamma<0$ then $\mathcal{S}_{2, \gamma}$ is not attained even if \eqref{bubblebis} are still solutions of the problem.\\

Here we are concerned with the quasilinear doubly critical case $1<p<N$ and $\gamma\in (0, \gamma_p)$. It is worth recalling that in
 \cite{AFP} the authors considered  minimization problem:
\begin{equation}\label{sgamma}\mathcal{S}_{p, \gamma}=\inf_{\substack{\varphi\in \mathcal D^{1, p}(\mathbb R^N) \\ \varphi\neq 0}}\frac{\displaystyle\intR \left(|\n \varphi|^p-\frac{\gamma}{|x|^p}\varphi^{p}\right)}{\left(\displaystyle\intR |\varphi|^{p^*}\right)^{\frac{p}{p^*}}}.\end{equation}
It follows that $0<\mathcal{S}_{p, \gamma}<\mathcal{S}_p$ where $\mathcal{S}_p$ is defined in \eqref{sp} and $\mathcal{S}_{p, \gamma}$ is attained by a function $u_0(x)$ which is not explicit. It has been proved in \cite{AFP}  that all minimizers of \eqref{sgamma} are radial. Also uniqueness up to scaling of the \emph{radial} solutions as well as  the asymptotic behavior are proved showing in particular that, given a \emph{radial} solution $u=u(r)$ to \eqref{pb}, then $$\lim_{r\to0} r^{\gamma_1} u(r)=C_1,\qquad \lim_{r\to+\infty}r^{\gamma_2}u(r)=C_2$$  and
\begin{equation*}
\lim_{r\to0} r^{\gamma_1+1} |u'(r)|=C_1\gamma_1,\qquad \lim_{r\to+\infty}r^{\gamma_2+1}|u'(r)|=C_2\gamma_2,
\end{equation*}
for some positive constants $C_1, C_2$. Here and hereafter  $\gamma_1, \gamma_2\in [0, +\infty)$, $\gamma_1<\gamma_2$ are defined as the two roots of the equation
\begin{equation}\label{radici}
\mu^{p-2}\left[(p-1)\mu^2-(N-p)\mu\right]+\gamma=0.
\end{equation}
We remark (for later use) that
\begin{equation*}\label{relgammaj}
0\le \gamma_1<\frac{N-p}{p}<\gamma_2\le \frac{N-p}{p-1}.
\end{equation*}
Note that when $p=2$ then $\gamma_{1}=\beta_-$ and $\gamma_2=\beta_+$. Instead, when $p\neq 2$ but $\gamma=0$ then $\gamma_1=0$ and $\gamma_2=\frac{N-p}{p-1}$.
Moreover in \cite{XiangJDE, Xiang2} the author extends the results on the asymptotic behavior proved for radial solutions in \cite{AFP} to all weak positive solutions of \eqref{pb}.\\

 We shall prove here that actually all positive solutions to \eqref{pb} are \emph{radially symmetric} thus allowing to deduce that the characterization of the solutions described here above do apply to all positive solutions. In particular, as a consequence of our result, we deduce uniqueness up to scaling of the positive solutions as well as  the their asymptotic behavior at the origin and at infinity.\\

 Our main result is the following:

\begin{theorem}\label{main}
Assume $\gamma\in (0, \gamma_p)$ and let $u$ be a positive solution to \eqref{pb}. Then $u$ is radial and radially decreasing with respect to the origin.
\end{theorem}

All the proofs of the classification results described above are based on the use of the {\it the moving plane method}. When $p\neq2$ this is completely not trivial because of the nonlinear degenerate nature of the operator. In our case, when trying to adapt the techniques developed in \cite{DR,DS1,S}, an obstruction occurs due to the homogeneity of the Hardy potential. In particular this fact is related to the nonlinear nature of the operator that also  obstructs the application of the techniques introduced in \cite{DMMS,Terracini}. In fact, to face this fact, we  exploit a different test function technique that, on the other hand, introduces  several difficulties as the reader shall see. Let us also stress that, for the absence of the Kelvin transformation, an analysis on the behaviour at infinity is needed. We will in fact exploit the results in \cite{XiangJDE, Xiang2} and in particular our Theorem \ref{stimagrad}.

\subsection{Notations}
Throughout the paper, we denote by $\Omega^c$ the complement of a domain $\Omega\subset \R$ in $\R$, by
$$C^k_0(\mathbb R^N)=\left\{u\in C^k(\mathbb R^N)\,: \ u(x)\to 0 \  \mbox{as} \ |x|\to+\infty\right\},$$
and by $B_R(x_0)$ the ball of radius $R$ centered at $x_0\in \R$.
\\ Moreover $\chi_{\Omega}$ is the characteristic function of the set $\Omega$, $(v-w)^+:=\max\{v-w, 0\}$ and $(v-w)^-:=\min\{v-w, 0\}$.
\\ Finally we underline that we will denote by $C, C_i, c_i$ several constants whose value
may change from line to line and, sometimes, on the same line. However these values will be not relevant in the proofs.

We remark that the potential $|x|^{-p}$ is related to the Hardy-Sobolev inequality. More precisely, for all $u\in \mathcal D^{1, p}(\mathbb R^N)$, one has
\begin{equation}\label{hs}
\intR \frac{|u|^p}{|x|^p}\le \frac{1}{\gamma_p}\intR |\nabla u|^p,\end{equation}
where $\gamma_p^{-1}$ is optimal and never achieved.
\\ As a consequence of a Pohozaev type identity, one can see that  problem \eqref{pb} does not have non-trivial solutions in any bounded starshaped domain with respect to the origin (Lemma 3.7 in \cite{GAP}).  \\

\section{Preliminaries and known technical results}\label{sezprel}
In this section we first recall useful results such as the strong comparison principle, a weighted Hardy-Sobolev inequality and decay estimates. \\

Let us start the discussion on the \emph{strong comparison principles} recalling the following
\begin{theorem}[Theorem $1.4$ of \cite{DS}]\label{DS}
 Let $u, v\in C^1(\bar \Omega)$ where $\Omega$ is a
bounded smooth domain of $\R$ with $\frac{2N+2}{N+2} <p< 2$ or $p > 2$. Suppose that
either $u$ or $v$ is a weak solution of
\begin{equation}\label{1.1}
\begin{cases}
-\Delta_p u =f(x,u) & \rm{in}\ \Omega,\\
u>0 & \rm{in}\ \Omega,\\
u=0 & \rm{on}\ \partial\Omega.\\
\end{cases}
\end{equation}
 with $f: \overline\Omega\times [0,\infty)\to \mathbb R$ is a continuous function which is positive and of class $C^1$ in  $\Omega\times(0, \infty)$. Assume that
$$\displaystyle -\Delta_p u +\Lambda u \le -\Delta_p v +\Lambda v\quad \mbox{and}\quad u\le v\quad \mbox{in}\,\, \Omega,$$

where $\Lambda\in \mathbb R$. Then $u\equiv v$ in $\Omega$ unless $
u < v$ in $\O$.\\
\end{theorem}
Actually the assumption that $u$ or $v$ fulfil the zero Dirichlet boundary datum can be removed and local versions of Theorem \ref{DS} are available, see \cite{Sci1,Sci2}. On the contrary there are no results removing the assumption $p>\frac{2N+2}{N+2}$. Therefore  in some cases we could prefer to exploit also the following result:

\begin{theorem}[Theorem $1.4$ of \cite{D}] \label{strong}
	Suppose $\Omega$ is a domain in $\mathbb R^N$ and let $u,v\in C^1(\Omega)$ weakly satisfy
	\begin{equation*}\label{D1}
	-\Delta_p u +\Lambda u \le -\Delta_p v +\Lambda v \quad \mbox{and}\quad u\le v\quad \mbox{in}\,\, \Omega,
	\end{equation*}
	$1<p<\infty$ and denote by $Z_v^u:=\left\{x\in\Omega\,:\, \nabla u(x)=\nabla v(x)=0\right\}$. Then if there exists $x_0\in \Omega\setminus Z_v^u$ with $u(x_0)=v(x_0)$, then $u\equiv v$ in the connected component of $\Omega\setminus Z_v^u$ containing $x_0$. The same result holds if more generally
	$$-\Delta_p u - f( u) \le -\Delta_p v - f(v)\quad \mbox{and}\quad u\le v\quad \mbox{in}\,\, \Omega,$$
	with $f:\mathbb R \to\mathbb R$ locally Lipschitz continuous.
\end{theorem}
In the spirit of the moving plane procedure we shall exploit the \emph{strong comparison principle} together with the \emph{weak comparison principle}
(that actually will be included in the proofs and we refer the readers to \cite{DS1}) and improved  Hardy inequalities proved in \cite{NS}.
For convenience we summarize the following
\begin{theorem}[Proposition 1.1 of \cite{NS}]\label{hardyns}
Let $r\ge 1$, $\tau>0$, $\alpha, \gamma\in \mathbb R$ such that
	\begin{equation*}\label{rel}
	\frac{1}{\tau}+\frac{\gamma}{N}=\frac{1}{r}+\frac{\alpha-1}{N},
	\end{equation*}
	and with $$0\le \alpha-\gamma\le 1.$$
	Let $u\in C^1_0(\mathbb{R}^N\setminus\{0\})$ and let $\frac 1 \tau +\frac \gamma N<0$ then it holds
		$$\left(\int_{ \mathbb R^N }|x|^{\gamma\tau}|u|^\tau\right)^{\frac 1 \tau}\le C \left(\int_{\mathbb{R}^N} |x|^{r\alpha}|\nabla u|^r\right)^{\frac{1}{r}} $$
	where $C$ is a positive constant independent of $u$.
\end{theorem}
\begin{remark}\label{remngu}
In Theorem \ref{hardyns} it is assumed that $u\in C^1_0(\mathbb{R}^N\setminus\{0\})$. Actually it is clear from the proof, and via density arguments, that the same result applies if $u$ is defined in exterior domains and has the right decay properties at infinity.
\end{remark}
To exploit Theorem \ref{hardyns}
for weak positive solutions to problem \eqref{pb} we need to know the
asymptotic behavior of the solution at infinity. Let us start recalling some results from
 \cite{XiangJDE, Xiang2}.
\begin{theorem}\label{xiang}
Let $u\in \mathcal D^{1,p}(\R)$ be a weak positive solution to equation \eqref{pb}. Then there exist positive constants $C, c$ depending on $N, p, \gamma$ and the solution $u$ such that
\begin{equation}\label{stimaviczero}
c|x|^{-\gamma_1}\le u(x) \le C|x|^{-\gamma_1}\qquad \mbox{for}\,\, |x|< R_0,
\end{equation}
and
\begin{equation}\label{stimalonzero}
c|x|^{-\gamma_2}\le u(x)\le C|x|^{-\gamma_2}\qquad \mbox{for}\,\, |x|> R_1.
\end{equation}
Moreover
\begin{equation}\label{gradaltoviczero}
|\n u(x)|\le c|x|^{-(\gamma_1+1)}\quad \ \ \ \ \ \ \ \ \mbox{for}\,\, |x|< R_0,
\end{equation}
and
\begin{equation}\label{gradaltolonzero}
|\n u(x)|\le c |x|^{-(\gamma_2+1)}\quad \ \ \ \ \ \ \ \ \mbox{for}\,\, |x|> R_1.
\end{equation}
Here $\gamma_1,\gamma_2$ are roots of \eqref{radici} and such that
\begin{equation*}\label{relgammajbis}
0\le \gamma_1<\frac{N-p}{p}<\gamma_2\le \frac{N-p}{p-1},
\end{equation*}
while $0< R_0<1< R_1$ are constants depending on $N, p, \gamma$ and the solution $u$.
\end{theorem}
Finally, we recall the following regularity result for solutions of \eqref{pb}.
\begin{theorem}[\cite{AFP,Di,T}]\label{reg}
Let $u$ be any solution of \eqref{pb}, then $u\in C^{1, \alpha}_{loc}(\mathbb R^N\setminus\{0\})$ with $0<\alpha<1$.
\end{theorem}

\section{Asymptotic estimates}
Here we shall prove some new gradient estimates  that we will use in the next section in order to apply the moving plane method. 
The moving plane procedure is strongly related to the use of suitable comparison principles. When the domain
is the whole space,  considering  problems  with a  source term involving the Hardy potential, weak comparison principles are naturally related to the use of Hardy type inequalities that involves the classical radial weights. Since our problem has a natural associated weight $|\nabla u|^{p-2}$, we will need to 
relate the weight $|\nabla u|^{p-2}$ with the weights appearing in  Theorem \ref{hardyns}. 
 To do this, especially  for the hardest case $p>2$, a further information is required, namely  estimates from below on the modulus of the gradient of the solution. This is what we prove in this section starting from the following:
\begin{lemma}
Let $u, v$ be positive and $C^1$-functions in a neighbourhood of some point $x_0\in \mathbb{R}^N$. Then it holds
\begin{equation}\label{stimamag2}
\begin{aligned}
&|\nabla u|^{p-2}\nabla u \cdot \nabla\left(u-\frac{v^p}{u^p} u\right)+|\nabla v|^{p-2}\nabla v \cdot \nabla\left(v-\frac{u^p}{v^p} v\right)
\\
& \ge C_p \min\{v^p,u^p\} \left(|\nabla\log u|+|\nabla \log v|\right)^{p-2}|\nabla\log u-\nabla\log v|^2,
\end{aligned}
\end{equation}
near $x_0$ for some constant $C_p$  depending only on $p$.
\end{lemma}
\begin{proof}
The estimate \eqref{stimamag2} for $1<p<2$ can be found in Lemma $3.1$ of \cite{XiangJDE}. Then we just need to prove \eqref{stimamag2} for $p\ge2$. \\
By making some simple computations we find that
\begin{equation}\label{stima0}
\begin{aligned}
T:&=|\nabla u|^{p-2}\nabla u \cdot \nabla\left(u-\frac{v^p}{u^p} u\right)+|\nabla v|^{p-2}\nabla v \cdot \nabla\left(v-\frac{u^p}{v^p} v\right) \\
&= |\nabla u|^p +|\nabla v|^p-v^p\underbrace{\left(|\nabla\log u|^p+p|\nabla\log u|^{p-2}\nabla\log u\cdot\(\nabla\log v-\nabla\log u\) \right)}_{(I)}\\
&-u^p\underbrace{\left(|\nabla\log v|^p+p|\nabla\log v|^{p-2}\nabla\log v\cdot\(\nabla\log u-\nabla\log v\) \right)}_{(II)}.
\end{aligned}
\end{equation}
Now let $f(t)=|a+t(b-a)|^p$ for $a, b\in \R$ then one has
$$f(1)=f(0)+f'(0)+\int_0^1 (1-t)f''(t),$$
which gives (recall that  $p\ge2$ ) 
\begin{equation}
\label{stima1}
\begin{aligned}
|b|^p&=|a|^p+p|a|^{p-2}a\cdot(b-a)\\
&+p(p-2)\int_0^1(1-t)|a+t(b-a)|^{p-4}\((a+t(b-a))\cdot(b-a)\)^2\, dt\\
&+p\int_0^1(1-t)|a+t(b-a)|^{p-2}|b-a|^2\,dt\\
&\ge |a|^p+p|a|^{p-2}a\cdot(b-a)+\int_0^1(1-t)p|a+t(b-a)|^{p-2}|b-a|^2\, dt.
\end{aligned}
\end{equation}
We apply \eqref{stima1} to (I) with $a=\nabla\log u$ and $b=\nabla\log v$ and to (II) with $a=\nabla\log v$ and $b=\nabla\log u$. Hence we get

\begin{equation}\label{T}
\begin{aligned}
T&\ge v^p\left[\int_0^1(1-t)p\left|\nabla\log u +t(\nabla \log v-\nabla\log u)\right|^{p-2}|\nabla\log u-\nabla\log v|^2\, dt\right]\\
&+ u^p\left[\int_0^1(1-t)p\left|\nabla\log v +t(\nabla \log u-\nabla\log v)\right|^{p-2}|\nabla\log u-\nabla\log v|^2\, dt\right] \\
&\ge \frac 34 p v^p|\nabla\log u-\nabla\log v|^2\left[\int_0^{\frac 14}\left|\nabla\log u +t(\nabla \log v-\nabla\log u)\right|^{p-2}\, dt\right] \\
&+\frac 34p u^p|\nabla\log u-\nabla\log v|^2\left[\int_0^{\frac 14}\left|\nabla\log v +t(\nabla \log u-\nabla\log v)\right|^{p-2}\, dt\right]. \\
\end{aligned}\end{equation}
Now suppose that $|\n \log u| \ge |\n \log v|$. In order to estimate the first term on the right hand side of \eqref{T} we distinguish two cases.
\\ First of all let $|\nabla\log v -\nabla\log u|\le \frac 12 |\nabla \log u|$ then (recall $0<t<1$)
$$\begin{aligned}|\n\log u+t(\n\log v-\n\log u)|&\ge |\n\log u|-|\n \log v-\n\log u|\\ &\ge \frac 12 |\n \log u|\ge \frac 14 \(|\n \log u|+|\n \log v|\),\end{aligned}$$
namely
$$|\n\log u+t(\n\log v-\n\log u)|^{p-2}\ge \(\frac 14\)^{p-2}\(|\n \log u|+|\n \log v|\)^{p-2}.$$
Otherwise if $|\nabla\log v -\nabla\log u|> \frac 12 |\nabla \log u|$ then we let $$t_0:=\frac{|\n\log u|}{|\n\log v-\n\log u|}\in (0, 2).$$ Hence
$$\begin{aligned}|\n\log u+t(\n\log v-\n\log u)|&\ge \left||\n \log u|-t|\n\log u-\n\log v|\right|\\
&=\left|t_0|\n\log u-\n\log v|-t|\n\log u-\n\log v|\right|\\
&=|t_0-t||\n\log u-\n\log v|\ge \frac 12 |t_0-t||\n\log u|\\
&\ge \frac 14 |t_0-t|\(|\n\log u|+|\n\log v|\),
\end{aligned}$$
since we are assuming that $|\n \log u| \ge |\n \log v|$. Therefore
$$|\n\log u+t(\n\log v-\n\log u)|^{p-2}\ge \(\frac 14\)^{p-2} |t_0-t|^{p-2}\(|\n\log u|+|\n\log v|\)^{p-2}.$$
Then, observing that $\int_0^{\frac 14}|t_0-t|^{p-2} \ge C_p$, one has
$$\begin{aligned} &\hskip-1.5cm \frac 34 pv^p|\nabla\log u-\nabla\log v|^2\left[\int_0^{\frac 14}\left|\nabla\log u +t(\nabla \log v-\nabla\log u)\right|^{p-2}\, dt\right]\\ & \ge C_p v^p\(|\n \log u|+|\n \log v|\)^{p-2}|\nabla\log u-\nabla\log v|^2. \end{aligned}$$
In the case $|\n \log u| \le |\n \log v|$, arguing in the same way, we deduce that
$$\begin{aligned} &\hskip-1.5cm\frac 34 p u^p|\nabla\log u-\nabla\log v|^2\left[\int_0^{\frac 14}\left|\nabla\log v +t(\nabla \log u-\nabla\log v)\right|^{p-2}\, dt\right]\\ & \ge C_p u^p\(|\n \log u|+|\n \log v|\)^{p-2}|\nabla\log u-\nabla\log v|^2, \end{aligned}$$
which concludes the proof.
\end{proof}
As we have already observed, a key tool in our proofs is the moving plane technique. To exploit it we need the following notations. We will study the symmetry of the solutions in the $\nu-$ direction for any $\nu\in S^{N-1}$ (i.e. $|\nu|=1$).
Since the problem is invariant up to rotations we fix $\nu=e_1$ and we let
$$T_\l=\left\{x\in \R\,:\, x_1=\l\right\},$$ $$\Sigma_\l=\left\{x\in \R\,:\, x_1<\l\right\},$$ $$x_\l= R_\l(x)=(2\l-x_1, x')\in \mathbb R\times \mathbb R^{N-1},$$ $$u_\l(x)=u(x_\l).$$ Now we state a result that will be used afterwards.
\begin{theorem}\label{pblim}
	Let $1<p<N$ and let $v \in C^{1,\alpha}_{loc}(\R \setminus \{0\})$ with $0<\alpha<1$ be a positive solution to
	\begin{equation}\label{eqV}-\Delta_p v - \frac{\gamma}{|x|^p}v^{p-1}=0\quad \mbox{in}\, \ \ \R\setminus \{0\},\end{equation}
	such that
	\begin{equation}\label{dec}\displaystyle \lim_{|x|\to 0} v(x)=\infty.\end{equation}
Then, if $v$ fulfils \eqref{stimalonzero}, it follows that
 $v$ is a radial (strict) decreasing function.
\end{theorem}
\begin{proof}
	First of all we need to prove that $v$ is a radial non-increasing function by applying the moving plane technique. We fix a direction $\nu =e_1$ and, for $\lambda<0$, we take as test function $\varphi_{1, \l}= v^{1-p}(v^p-v_\l^p)^+\chi_{\Sigma_\l}$ and $\varphi_{2, \l}=v_\l^{1-p}(v^p-v_\l^p)^+\chi_{\Sigma_\l}$ in the weak formulation solved, respectively, by $v$ and $v_\l$. We note that $v_\lambda$ solves
\begin{equation}\label{vlambda} -\Delta_p v_\l -\frac{\gamma}{|x_\l|^p}v_\l^{p-1}=0.\end{equation} We also remark that, by using \eqref{dec},  $${\rm supp}(\varphi_{j, \lambda})\subset\subset \Sigma_\lambda\setminus\{0_\lambda\}\qquad j=1, 2.$$
It is easy to verify that $\varphi_{1, \l}\,,\,\varphi_{2, \l}\in\mathcal D^{1, p}(\mathbb R^N)$ (here we also exploit \eqref{stimalonzero}).
Furthermore, since  $\varphi_{1, \l}\,,\,\varphi_{2, \l}$ have compact
support far from the singularities,  we can use the weak formulations of
 \eqref{eqV}, \eqref{vlambda}  and, taking the difference, we deduce that
	\begin{equation}\label{stimarad}
	\begin{aligned}
	&\int_{\Sigma_\l}|\nabla v|^{p-2}\nabla v\cdot\nabla\varphi_{1, \l}-|\nabla v_\l|^{p-2}\nabla v_\l\cdot\nabla\varphi_{2, \l}	
	\\
	& +\gamma\int_{\Sigma_\l}\left(-\frac{1}{|x|^p}+\frac{1}{|x_\l|^p}\right)(v^p-v_\l^p)^+=0,
	\end{aligned}
	\end{equation}
	and, since $|x|>|x_\l|$ in $\Sigma_\l$, the second term on the left hand side of \eqref{stimarad} is nonnegative. Then, exploiting \eqref{stimamag2}, it follows that
	\begin{equation*}
		C_p \int_{\Sigma_\l\cap\{v\ge v_\l\}} v_\l^p\left(|\nabla\log v|+|\nabla \log v_\l|\right)^{p-2}|\nabla(\log v-\log v_\l)|^2 = 0
	\end{equation*}
	which implies that $\log v -\log v_\l$ is constant $\Sigma_\l\cap\{v\ge v_\l\}$ and since $\log v -\log v_\l =0$ on $T_\l$ we have $v\le v_\l$ on $\Sigma_\l$ for any $\lambda<0$. We repeat the same argument in the $-e_1$ direction deducing that $v$ is symmetric with respect to the $e_1$-direction. This procedure can be clearly performed in any direction $\nu \in \mathcal{S}^{N-1}$ whence one gets the radial monotone nonincreasing behavior of $v$.\\

A simple application of the Hopf Lemma (that can be applied since the level sets are spheres) shows now that $v$ has no critical points and in particular the radial derivative is strictly negative.
\end{proof}

Next we provide the corresponding lower bound for the decay rate of $|\n u|$ of Theorem \ref{xiang}.
\begin{theorem}\label{stimagrad}
	Let $1<p<N$ and let $u$ be a solution of \eqref{pb}. Then there exists $ R_2>0$ and a constant $\bar C>0$ such that
	\begin{equation}\label{grad}
	|\nabla u(x)|\ge \frac{\bar C}{|x|^{\gamma_2+1}} \quad \mbox{for}\,\, |x|> R_2.\end{equation}
\end{theorem}
\begin{proof}
Once that Theorem \ref{pblim} is in force we can carry out the proof borrowing
some ideas from Theorem 2.2 of \cite{S} . We sketch it for the sake of completeness.\\
By contradiction let us assume that there exist sequences of radii $R_n$ and points $x_n$ with $R_n\to+\infty$  as $n\to +\infty$ and $|x_n| =R_n$, such that
\begin{equation}\label{grad1}
|\n u (x_n)|\le \frac{\theta_n}{|R_n|^{\gamma_2+1}},
\end{equation}
with $\theta_n\to 0$ as $n\to+\infty$.
Without loss of generality we suppose $R_n>1$ for any $n$ and we set $w_{R_n}(x):=R_n^{\gamma_2} u(R_n x)$. One can observe that for fixed $0<a<A$ then $||w_{R_n}||_{L^\infty(B_A\setminus B_a)}$ is bounded with respect to $n$. Otherwise if $|x|> \frac{R_1}{R_n}$ one deduces by Theorem \ref{xiang} that
$$\frac{\bar c}{A^{\gamma_2}}\le w_{R_n}(x)\le \frac{\bar C}{a^{\gamma_2}},$$ and that
\begin{equation}
\label{grad2}
\begin{cases}\vspace{0.2cm}
\displaystyle w_{R_n}(x)\le \frac{\bar C}{A^{\gamma_2}}\,\, &x\in  \partial B_A, \\
\displaystyle w_{R_n}(x)\ge \frac{\bar c}{a^{\gamma_2}}\,\, &x\in \partial B_a.
\end{cases}
\end{equation}
Therefore, the above bound in $L^\infty(B_A\setminus B_a)$ implies that $w_{R_n}$ is also uniformly bounded in $C^{1,\alpha}(K)$ with $0<\alpha<1$ for any compact set $K\subset B_A\setminus B_a$. Finally, since $a>0$, without loss of generality we suppose that the $C^{1,\alpha}$ estimates hold in the closure of $B_A\setminus B_a$. Hence, for $x\in B_A\setminus B_a$ and up to subsequences, one gets that $w_{R_n}(x) \longrightarrow w_{a, A}(x)$ in $C^{1, \alpha'}$ for $0<\alpha'<\alpha$. We also underline that $w_{a, A}(x)$ satisfies \eqref{grad2}.
Furthermore, since
$$-\Delta_p w_{R_n} -\frac{\gamma}{|x|^p}w_{R_n}^{p-1}=\frac{w_{R_n}^{p^*-1}}{R_n^{(p^*-p)\gamma_2-p}}\quad\mbox{in}\, \ \ \R,$$
then
\begin{equation}\label{eqwaa}
	-\Delta_p w_{a, A} -\frac{\gamma }{|x|^p} w_{a, A}^{p-1}=0\qquad \mbox{in}\,\, B_A\setminus \overline{B_a}.
\end{equation}
Now, for $j\in \mathbb N$, one can take $a_j =\frac 1 j$ and $A_j = j$ and reasoning as above one constructs $w_{a_j,A_j}$.
Then, for $j\to\infty$, a diagonal argument implies the existence of a limiting profile $w_\infty$ such that $w_\infty\equiv w_{a_j,A_j}$ in $B_{A_j}\setminus B_{a_j}$.
In particular from \eqref{eqwaa} read for $w_{a_j, A_j}$ one has
$$-\Delta_p w_\infty -\frac{\gamma }{|x|^p}w_\infty^{p-1} =0\qquad \mbox{in}\,\, \mathbb R^N\setminus \{0\}.$$
From \eqref{grad2} with $a =a_j$ and $A =A_j$, one gets that the limiting profile $w_\infty$ is such that
$$\lim_{|x|\to+\infty}w_\infty(x) =0\qquad \mbox{and}\qquad \lim
_{|x|\to0} w_\infty (x) = +\infty$$
and it satisfies \eqref{stimalonzero}. Therefore
 Theorem \ref{pblim} can be applied providing that $w_\infty$ is radial with negative radial derivative.\\

To conclude let now $x_n$ be as in \eqref{grad1} and set
$y_n =\frac{ x_n}{R_n}.$ Then, by \eqref{grad1}, it follows that $|\n w_{R_n}(y_n)|$ tends to zero as $n\to+\infty$. Up to subsequences, since $|y_n| =1$, we have that $y_n \to \bar{y}\in \partial B_1$. Consequently, by the uniform convergence of the gradients one has that $\n w_\infty(\bar y) = 0$,  which is in contradiction with the definition of $w_\infty$, since, by Theorem \ref{pblim}, this cannot happen. \bk
\end{proof}
%
%
%

\section{Proof of the symmetry result}\label{princ}
We are now able to prove Theorem \ref{main}. First of all we underline that it is easy to see that $u_\l$ solves
\begin{equation}\label{pbriflessa}
-\Delta_p u_\l-\frac{\gamma}{|x_\l|^p}u_\l^{p-1}=u_\l^{p^*-1}\quad \mbox{in} \ \R.\end{equation} In what follows we set
$$\L^-=\left\{\l<0\,:\, u\le u_\mu\,\, \mbox{in}\,\, \Sigma_\mu, \,\, \forall\, \mu\le \l\right\},\quad \L^+=\left\{\l>0\,:\, u\ge u_\mu\,\, \mbox{in}\,\, \Sigma_\mu,\,\, \forall\, \mu\le \l\right\}.$$
If $\L^-\neq \emptyset$ and $\L^+\neq \emptyset$ we denote by $\l^-_0:=\sup\L^-$ and by $\l^+_0:=\inf\L^+$.\\

Roughly speaking, the moving plane method consists of two main steps: first in reflecting the domain about a fixed hyperplane and proving that the value the solution at each reflected point is larger than the value at the point itself and secondly in moving the hyperplane to a critical position; finally the solution results to be symmetric with respect to this limit hyperplane.
\begin{proof}[Proof of Theorem \ref{main}]
	We prove the result by analizing, sometimes in different ways, the case $1<p<2$ and the case $p>2$. For $p=2$ we refer to \cite{Terracini}. We divide the proof in two steps. \\ \\
{\bf Step 1:} 	$\L^-\neq \emptyset$ and $\L^+\neq \emptyset$.
\\ \\ We only prove $\L^-\neq \emptyset$, which is the existence of $\l<0$ with $|\l|$ sufficiently large such that  $u\le u_\mu$ in $\Sigma_\mu$ for every $\mu\le \l$. The proof of the fact that $\L^+\neq \emptyset$ is analogous and, at the end of the step, we outline the main changes in the proof in order to conclude it.\\
For the entire proof we denote by $R_0$, $R_1$ and $R_2$ the radii given by \eqref{stimaviczero}, \eqref{stimalonzero} and \eqref{grad} and we firstly observe that for $|\overline{\l}|>\max(R_1,R_2)$ one has, by \eqref{stimaviczero} and \eqref{stimalonzero}, that there exists $\tilde{R}_0:=\tilde{R}_0(\overline{\l})$ such that $\tilde{R}_0<R_0$, $B_{\tilde{R}_0}(0_{\overline{\l}})\subset \Sigma_{\overline{\l}}$ and
\begin{equation}\label{stimartilde}
	\displaystyle \sup_{x\in B_{\tilde{R}_0}(0_{\overline{\l}})} u(x) < \displaystyle \inf_{x\in B_{\tilde{R}_0}(0_{\overline{\l}})}u_{\overline{\l}}(x).
\end{equation}
Therefore,  exploiting also \eqref{stimalonzero}, we deduce that
$$\displaystyle \sup_{x\in B_{\tilde{R}_0}(0_\l)} u(x) \leq \displaystyle \inf_{x\in B_{\tilde{R}_0}(0_\l)}u_{\l}(x),$$
which gives that $u<u_\l$ in $B_{\tilde{R_0}}(0_\l)\subset\Sl$ for every $\l\le\overline{\l}$ and with $\tilde{R_0}$ independent of $\l$.
Moreover we also denote by $\eta\in C^\infty_0(B_{2R}(0))$ a cut-off function such that $0\le \eta\le 1$, $\eta\equiv 1$ on $B_R(0)$ and $|\nabla \eta|\le \frac{2}{R}$.
\\ \\ In what follows we employ the following notation: $\Sigma_\l'=\Sigma_\l\setminus B_{\tilde{R}_0}(0_\l)$ and $\hat B_{\rho}:=B_\rho(0)\cap \Sigma_\l'$ for $\rho>0$.\\
\\ If $\alpha>\max\{2, p\}$ and $\l \le \overline{\l}$, we consider
\begin{equation}\label{test}
\varphi_{1, \l}=\eta^\alpha u^{1-p}(u^p-u_\l^p)^+\chi_{\Sigma_\l},\qquad \varphi_{2, \l}=\eta^\alpha u_\l^{1-p}(u^p-u_\l^p)^+\chi_{\Sigma_\l}.
\end{equation}
We remark that ${\rm supp}(\varphi_{j, \l}) \subset \hat B_{2R}$ for $j=1, 2$.
 Then we take $\varphi_{1, \l}$ as a test function in \eqref{pb}, $\varphi_ {2, \l}$ in \eqref{pbriflessa} and we subtract.
Hence, denoting by $\psi_\l:=(u^p-u_\l^p)^+$ and by $\varphi_\l:=(u-u_\l)^+$ one gets
\begin{equation}\label{sim1}
\begin{aligned}
&\int_{\hat B_{2R}}\left(|\nabla u|^{p-2}\nabla u\cdot\nabla\varphi_{1, \l}-|\nabla u_\l|^{p-2}\nabla u_\l\cdot\nabla\varphi_{2, \l}\right)
+\gamma\int_{\hat B_{2R}}\left(-\frac{1}{|x|^p} + \frac{1}{|x_\l|^p}\right)\eta^\alpha\psi_\l\\
&=\int_{\hat B_{2R}}(u^{p^*-p}-u_\l^{p^*-p})\eta^\alpha\psi_\l,
\end{aligned}
\end{equation}
and, since $|x|\ge |x_\l|$ in $\Sl$, one has that the second term on the left hand side of \eqref{sim1} is nonnegative. Hence
\begin{equation}\label{inizio}
\begin{aligned}
&\underbrace{\int_{\hat B_{2R}}\eta^\alpha \left(|\nabla u|^{p-2}\nabla u\cdot\nabla(u^{1-p}\psi_\l)-|\nabla u_\l|^{p-2}\nabla u_\l\cdot\nabla(u_\l^{1-p}\psi_\l)\right)}_{I_1}\\
&\le\underbrace{-\alpha\int_{\hat B_{2R}}\eta^{\alpha-1}u^{1-p}\psi_\l |\nabla u|^{p-2}\nabla u\cdot \nabla \eta}_{I_2} +\underbrace{\alpha\int_{\hat B_{2R}}\eta^{\alpha-1}u_\l^{1-p}\psi_\l|\nabla u_\l|^{p-2}\nabla u_\l \cdot \nabla \eta}_{I_3} \\
&+\underbrace{\int_{\hat B_{2R}}(u^{p^*-p}-u_\l^{p^*-p})\eta^\alpha\psi_\l}_{I_4}.
\end{aligned}
\end{equation}
We start by estimating $I_1$. By using \eqref{stimamag2} it yields that for $p>2$ one has
\begin{equation}\label{I1}
\begin{aligned}
I_1&\ge C_p \int_{\hat B_{2R}\cap \{u\ge u_\l\}} \eta^\alpha u_\l^p\(|\n \log u|+|\n\log u_\l|\)^{p-2}|\n\log u-\n\log u_\l|^2
\\
&\ge C_p\int_{\hat B_{2R}\cap \{u\ge u_\l\}}\eta^\alpha\(\frac{u_\l}{u}\)^p u^2\(|\n u|+|\n u_\l|\)^{p-2}|\n\log u -\n\log u_\l|^2
\\
&\ge c_1 \int_{\hat B_{2R}\cap \{u\ge u_\l\}}\eta^\alpha u^2\(|\n u|+|\n u_\l|\)^{p-2}|\n\log u -\n\log u_\l|^2,
\end{aligned}
\end{equation}
while for $1<p<2$ we obtain
\begin{equation}\label{I1bis}
\begin{aligned}
I_1&\ge C_p \int_{\hat B_{2R}\cap \{u\ge u_\l\}} \eta^\alpha  u_\l^p\frac{|\n\log u-\n\log u_\l|^2}{(|\n \log u|+|\n \log u_\l|)^{2-p}}
\\
&\ge C_p\int_{\hat B_{2R}\cap \{u\ge u_\l\}}\eta^\alpha u_\l^2 \frac{|\n\log u -\n\log u_\l|^2}{\(|\n u|+|\n u_\l|\)^{2-p}}.
\end{aligned}
\end{equation}
We remark that in \eqref{I1} we used that
\begin{equation}\label{stimaulu}
\frac{u_\l}{u}\ge \tilde c \quad \mbox{in}\ \Sigma_\l,
\end{equation}
and $c_1:=C_p\tilde{c}^p$.
Indeed if $x \in \Sl\setminus B_{R_1}(0_\lambda)$ then from \eqref{stimaviczero} and \eqref{stimalonzero} one has (recall that $|x|\ge |x_\l|$)
\begin{equation*}
	\frac{u_\l}{u}\ge \tilde c_1\frac{|x|^{\gamma_2}}{|x_\l|^{\gamma_2}}\ge \tilde c_1.
\end{equation*}
Otherwise if $x \in  \Sl\cap B_{R_1}(0_\lambda)$ then
\begin{equation*}
\frac{u_\l}{u}\ge \tilde c_1|\overline{\l}|^{\gamma_2} \inf_{x\in B_{R_1}(0)} u(x) \ge \tilde c_2,
\end{equation*}
and we set $\tilde c=\min(\tilde c_1,\tilde c_2)$.
Now it follows from \eqref{stimalonzero} and \eqref{gradaltolonzero} that
\begin{equation}\label{I2finale}
\begin{aligned} I_2 &\le \alpha \int_{\hat B_{2R}\cap \{u\ge u_\l\}}\eta^{\alpha-1}u\left(1-\(\frac{u_\l}{u}\)^p\right)|\nabla u|^{p-1}|\n\eta|\\
&\le \frac{2\alpha}{R}\int_{\hat B_{2R}\setminus \hat B_R}u|\n u|^{p-1}  \le \frac{C}{R}\int_{\hat B_{2R}\setminus \hat B_{R}}\frac{1}{|x|^{(\gamma_2+1)(p-1)+\gamma_2}} \le \frac{C}{R^\beta}
\end{aligned}
\end{equation}
where, from here on, $\beta:=p\gamma_2+p-N$ which is strictly positive since $\gamma_2>\frac{N-p}{p}$. For $I_3$, using \eqref{stimalonzero} and \eqref{stimaulu}, we deduce that
\begin{equation}\label{I3finale}
\begin{aligned} I_3&\le C\int_{\hat B_{2R}\setminus \hat B_R} \alpha \eta^{\alpha-1}u_\l^{1-p}\(u^p-u_\l^p\)^+|\n u_\l|^{p-1}|\n \eta|
\displaystyle
\\
&\le \frac 2 R \int_{\hat B_{2R}\setminus \hat B_R\cap \{u\ge u_\l\}}u_\l\(\(\frac{u}{u_\l}\)^p-1\)|\n u_\l|^{p-1}
\\
&\le \frac 2 R \int_{\hat B_{2R}\setminus \hat B_R\cap \{u\ge u_\l\}}u_\l \(\frac{u}{u_\l}\)^p |\n u_\l|^{p-1}\le \frac C R  \int_{\hat B_{2R}\setminus \hat B_R\cap \{u\ge u_\l\}}u|\n u_\l|^{p-1}\\
&\le \frac C R  \(\int_{\R} |\n u_\l|^p\)^{\frac{p-1}{p}}\(\int_{\hat B_{2R}\setminus \hat B_R}u^p\)^{\frac 1 p}\le \frac C R \(\int_{\hat B_{2R}\setminus \hat B_R}\frac{1}{|x|^{\gamma_2p}}\)^{\frac 1 p}\le \frac {C}{ R^\frac{\beta}{p}}.
\end{aligned}
\end{equation}
For the term $I_4$ we first note that  (since $u\ge u_\l$)
$$\begin{aligned}I_4&= \int_{\hat B_{2R}}\left(\frac{u^{p^*-1}}{u^{p-1}}-\frac{u_\l^{p^*-1}}{u_\l^{p-1}}\right)\eta^\alpha \psi_\l \le  \int_{\hat B_{2R}}\frac{1}{u^{p-1}}\left(u^{p^*-1}-u_\l^{p^*-1}\right)\eta^\alpha \psi_\l,\end{aligned}$$
then applying twice the Lagrange Theorem and using \eqref{stimalonzero} one has that in case  $p^*\ge 2$
\begin{equation*}\label{I40}
I_4 \le c_p\int_{\hat B_{2R}}  u^{p^*-2}\eta^\alpha\varphi_\l^2\le c_p \int_{\hat B_{2R}} \frac{1}{|x|^{\gamma_2(p^*-2)}}\eta^\alpha\varphi_\l^2,
\end{equation*}
while for $1<p^*<2$ (recall \eqref{stimaulu})
\begin{equation*}\label{I42}
I_4  \le c_p\int_{\hat B_{2R}} \frac{\eta^\alpha\varphi_\l^2}{u_\l^{2-p^*}}=\int_{\hat B_{2R}}\left(\frac{u}{u_\l}\right)^{2-p^*} \frac{\eta^\alpha\varphi_\l^2}{u^{2-p^*}}\le c_p\int_{\hat B_{2R}} \frac{\eta^\alpha\varphi_\l^2}{u^{2-p^*}}\le c_p \int_{\hat B_{2R}} \frac{1}{|x|^{\gamma_2(p^*-2)}}\eta^\alpha\varphi_\l^2,
\end{equation*}
which gives  for any $p>1$
\begin{equation}\label{I4}
I_4 \le c_p \int_{\hat B_{2R}} \frac{1}{|x|^{\gamma_2(p^*-2)}}\eta^\alpha\varphi_\l^2.
\end{equation}
Let us now consider $f(t)=\log(a+t(b-a))$ where $a,b>0$ ($b\ge a$) then

$$\log b=\log a +(b-a)\int_0^1 \frac{1}{a+t(b-a)},$$ and since $t\in [0, 1]$ we get

\begin{equation}\label{usare} b-a = \frac{\log b-\log a}{\int_0^1 \frac{1}{a+t(b-a)}}\le b (\log b-\log a).\end{equation}

We use \eqref{usare} with $b=u$ and $a=u_\l$ and estimate the right hand side of \eqref{I4} (by using also \eqref{stimalonzero}) as
\begin{equation*}\label{i42}
\begin{aligned}
	I_4&\le C \int_{\hat B_{2R}\cap \{u\ge u_\l\}}\frac{1}{|x|^{\gamma_2(p^*-2)}}\eta^\alpha u^2\(\log u-\log u_\l\)^2
	\\
	&\le C \int_{\hat B_{2R}}\frac{1}{|x|^{\gamma_2 p^*}}\eta^\alpha \left((\log u-\log u_\l)^+\right)^2.\\
\end{aligned}
\end{equation*}
Moreover
\begin{equation}\label{I4bis}
\begin{aligned}
	I_4 &\le C \int_{\hat B_{2R}}\frac{1}{|x|^{\beta^*-2\alpha+2}}\left(\eta^{\frac \alpha 2}(\log u-\log u_\l)^+\right)^2
	\\
	&\le  \frac{C}{|\l|^{\beta^*}} \int_{\hat B_{2R}}|x|^{2\alpha-2}\left(\eta^{\frac \alpha 2}(\log u-\log u_\l)^+\right)^2,
\end{aligned}
\end{equation}
where
$$\beta^*:=\gamma_2(p^*-p)-p;\qquad 2\alpha:=-[(\gamma_2+1)(p-2)+2\gamma_2].$$
We underline that $\beta^*-2\alpha+2=\gamma_2p^*$ and that $\beta^*>0$ since $\gamma_2>\frac{N-p}{p}$.
For the right hand side of \eqref{I4bis} we can apply Theorem \ref{hardyns} where $r=2, \tau=2$ which implies that
$$\gamma:=\alpha-1=-\frac{(\gamma_2+1)p}{2}$$
and that
$$\frac 1 2 +\frac \gamma N =\frac{N-\gamma_2p-p}{2N}<0$$ since $\gamma_2>\frac{N-p}{p}$.  Hence we obtain
\begin{equation}\label{stimaI4}
I_4 \le \frac{C}{|\l|^{\beta^*}} \int_{\hat{B}_{2R}}|x|^{2\alpha}|\nabla (\eta^{\frac \alpha 2}(\log u -\log u_\l)^+)|^2,
\end{equation}
and now, in order to estimate the right hand side of \eqref{stimaI4}, we distinguish between the case $p>2$ and the case $1<p<2$.
From \eqref{stimaI4} and for $p>2$ we get
\begin{equation}\label{stimaI4pmag2}\begin{aligned}
I_4 &\le \frac{C}{|\l|^{\beta^*}} \int_{\hat{B}_{2R}\cap\{u\ge u_\l\}}\frac{1}{|x|^{(\gamma_2+1)(p-2)}}\eta^\alpha u^2|\nabla \log u - \n \log u_\l|^2
\\
&+\frac{C}{|\l|^{\beta^*}} \int_{\hat{B}_{2R}\cap\{u\ge u_\l\}}|x|^{2\alpha}\(\log u - \log u_\l\)^2|\n \eta|^2
\\
&\le \frac{C}{|\l|^{\beta^*}} \int_{\hat{B}_{2R}\cap\{u\ge u_\l\}}\eta^\alpha u^2|\n u|^{p-2}|\nabla \log u - \n \log u_\l|^2+ \frac{C}{|\l|^{\beta^*}R^2} \int_{\hat{B}_{2R}\setminus\hat{B}_{R}}|x|^{2\alpha}
\\
&\le \frac{C}{|\l|^{\beta^*}} \int_{\hat{B}_{2R}\cap\{u\ge u_\l\}}\eta^\alpha u^2\(|\n u|+|\n u_\l|\)^{p-2}|\nabla \log u - \n \log u_\l|^2+ \frac{C}{|\l|^{\beta^*}R^\beta}.
\end{aligned}\end{equation}

Then, by using the estimates \eqref{I1}, \eqref{I2finale},  \eqref{I3finale} and  \eqref{stimaI4pmag2} in \eqref{inizio}, we
$$\(c_1-\frac{C}{|\l|^{\beta^*}}\) \int_{\hat B_{2R}\cap \{u\ge u_\l\}}\eta^\alpha u^2\(|\n u|+|\n u_\l|\)^{p-2}|\n\log u -\n\log u_\l|^2\le \frac{C}{R^\frac{\beta}{p}}  +\frac{C}{|\l|^{\beta^*}R^\beta}+\frac{C}{R^\beta}.$$
For $|\l|$ sufficiently large, as $R$ goes to $+\infty$, we deduce that 
$$\begin{aligned}&\int_{\Sigma'_\l \cap \{u\ge u_\l\}}u^2\( |\nabla u|+|\n u_\l|\)^{p-2} |\n \log u- \n \log u_\l|^2\\ &=  \lim_{R\to+\infty}  \int_{\hat B_{R} \cap \{u\ge u_\l\}}u^2\( |\nabla u|+|\n u_\l|\)^{p-2}|\n \log u- \n \log u_\l|^2\le 0.\end{aligned}$$
Now we have to estimate the right hand side of \eqref{stimaI4} in the case $1<p<2$.\\
\\ We first remark that $2\alpha<0$ (for $N>2$) and, since $|x|\ge |x_\l|$, one has that $|x|^{2\alpha}\le |x_\l|^{2\alpha}$.
Then
\begin{equation}\label{stimaI4pmin2}\begin{aligned}
I_4 &\le \frac{C}{|\l|^{\beta^*}}\int_{\hat{B}_{2R}\cap\{u\ge u_\l\}}|x_\l|^{2\alpha}\eta^\alpha |\nabla \log u -\nabla \log u_\l|^2
\\
&+\frac{C}{|\l|^{\beta^*}}\int_{\hat{B}_{2R}\cap\{u\ge u_\l\}}|x|^{2\alpha}\(\log u - \log u_\l\)^2|\n \eta|^2.
\end{aligned}\end{equation}
Let $\bar R=\max\{R_1, R_2\}$ and let $A_{\bar R, \tilde R_0}=\overline{B_{\bar R}(0_\l)\setminus B_{\tilde R_0}(0_\l)}$. Then we get $$\hat B_{2R}= \hat A_{\bar R, \tilde R_0}\cup \(\hat B_{2R}\setminus \hat A_{\bar R, \tilde R_0}\).$$
Exploiting  \eqref{stimalonzero} we deduce that
\begin{equation}\label{stimaI4pmin2due}\begin{aligned}
&\int_{\hat B_{2R}\setminus \hat A_{\bar R, \tilde R_0}}|x_\l|^{2\alpha}\eta^\alpha |\nabla \log u - \n \log u_\l|^2
\\
&\le C \int_{\hat B_{2R}\setminus \hat A_{\bar R, \tilde R_0}}|x_\l|^{(\gamma_2+1)(2-p)}|x_\l|^{-2\gamma_2}\eta^\alpha |\nabla \log u - \n \log u_\l|^2
\\
&\le C \int_{\hat B_{2R}\setminus \hat A_{\bar R, \tilde R_0}}u_\l^2\eta^\alpha\frac{ |\nabla \log u - \n \log u_\l|^2}{\(|\n u|+|\n u_\l|\)^{2-p}}.
 \end{aligned}\end{equation}
In $ A_{\bar R, \tilde R_0}$ it holds that $|x_\l| \ge \tilde R_0 $ and, since we are far from $0_\l$, we also get that $|\n u_\l|$ is bounded. Let $L:=\displaystyle \inf_{B_{\bar R}(0)\setminus B_{\tilde R_0}(0)} u$. Hence we get (by using \eqref{stimaulu} and the fact that $\(|\n u|+|\n u_\l|\)^{2-p}\le C$ away from $0, 0_\l$)
\begin{equation}\label{stimaI4pmin2tre}
\begin{aligned}\int_{A_{\bar R, \tilde R_0}}|x_\l|^{2\alpha}\eta^\alpha |\nabla \log u - \n \log u_\l|^2 &\le C\tilde R_0^{2\alpha}\int_{A_{\bar R, \tilde R_0}}\eta^\alpha |\nabla \log u - \n \log u_\l|^2\\
&\le \frac{C\tilde R_0^{2\alpha}}{L^2}\int_{A_{\bar R, \tilde R_0}}u_\l^2\eta^\alpha \frac{ |\nabla \log u - \n \log u_\l|^2}{\(|\n u|+|\n u_\l|\)^{2-p}}\,. \end{aligned}
\end{equation}
Gathering \eqref{stimaI4pmin2due} and \eqref{stimaI4pmin2tre} in the first term of \eqref{stimaI4pmin2} and reasoning as in \eqref{stimaI4pmag2} for the second term of \eqref{stimaI4pmin2} one yields to 
\begin{equation}\label{stimaI4pmin2unoprimo}\begin{aligned}
&I_4 \le \frac{C}{|\l|^{\beta^*}} \int_{\hat{B}_{2R}\cap \{u\ge u_\l\}}u_\l^2 \eta^\alpha \frac{ |\nabla\log u - \n\log u_\l|^2}{\(|\n u|+|\n u_\l|\)^{2-p}} + \frac{C}{|\l|^{\beta^*}R^\beta}. \end{aligned}\end{equation}
Hence, by collecting \eqref{I1bis}, \eqref{I2finale}, \eqref{I3finale} and \eqref{stimaI4pmin2unoprimo} in \eqref{inizio}, we get
$$\(c_1-\frac{C}{|\l|^{\beta^*}}\) \int_{\hat B_{2R}\cap \{u\ge u_\l\}}\eta^\alpha u_\l^2\frac{ |\nabla \log u - \n \log u_\l|^2}{\(|\n u|+|\n u_\l|\)^{2-p}}\le \frac {C}{ R^\frac{\beta}{p}}+\frac {C}{ R^\beta} + \frac{C}{|\l|^{\beta^*}R^\beta} .$$
Once again we can choose $|\l|$ large enough so that, as $R$ goes to $+\infty$, it yields
$$\begin{aligned}&\int_{\Sigma'_\l \cap \{u\ge u_\l\}}u_\l ^2\frac{ |\nabla \log u - \n \log u_\l|^2}{\(|\n u|+|\n u_\l|\)^{2-p}}=\limsup_{R\to+\infty}\int_{\hat B_{R}}u_\l^2\frac{  |\nabla \log u - \n \log u_\l|^2}{\(|\n u|+|\n u_\l|\)^{2-p}}\le 0.\end{aligned}$$
Hence, in both cases, $\log u -\log u_\l$ is constant and since $\log u -\log u_\l =0$ on $T_\l$ then $\log u-\log u_\l =0$ on the set $\Sigma'_\l \cap \{u\ge u_\l\}$. Therefore we get $u\le u_\l$ on $\Sigma_\l$. Hence $\L^-\neq \emptyset$ and $\l^-_0$ exists and it is also finite.\\
In order to show that	$\L^+\neq \emptyset$ then we take as test functions $$\phi_{1, \l}=u^{1-p}(u^p-u_\l^p)^-\chi_{\Sigma_\l}, \ \phi_{2, \l}=u_\l^{1-p}(u^p-u_\l^p)^-\chi_{\Sigma_\l}$$
and, analogously to what already done, we are able to prove the claim so that there exists $\l_0^+$ which is also finite.
\\ \\
{\bf Step 2:} 	$\l_0^-=\l_0^+=0$. \\
We argue by contradiction assuming that $\l_0^-\neq 0$.  Arguing as in the proof of Step $1$  we will get the contradiction  proving that $u\le u_{\l_0^- +\varepsilon}$ in $\Sigma_{\l_0^-+\varepsilon}$ for all $0\leq\varepsilon \le \overline{\varepsilon}$ for some $\overline{\varepsilon}>0$.\\

In what follows we shall exploit the \emph{strong comparison principle}. To do this we start noticing that
from Step 1 and by continuity it holds that
$$u\le u_{\l_0^-}\quad \text{in}\quad\Sigma_{\l_0^-}\,.$$
By Theorem \ref{strong} we deduce that
$u\equiv u_{\l_0^-}$ or $u<u_{\l_0^-}$ in any connected component $\mathcal C$ of $\Sigma_{\l_0^-}\setminus{Z_u}$ ($Z_u=\{\nabla u=0\}$).
We will frequently use the fact that  $Z_u$ has zero  Lebesgue measure \cite{DS1}.\\

\noindent Assume first that $\Sigma_{\l_0^-}\setminus{Z_u}$ has only one connected component.  We observe that $u\equiv u_{\l_0^-}$ is not possible in this case since, by \eqref{stimaviczero}, there exists $B_{\tilde R_0}(0_{\l_0^-})$ where $u<u_{\l_0^-}$; this means that $u<u_{\l_0^-}$ in $\Sigma_{\l_0^-}\setminus Z_u$.\\

Assume now that there are at least two connected components of $\Sigma_{\l_0^-}\setminus{Z_u}$. Our Theorem \ref{stimagrad} implies that $Z_u$ is bounded so that only one component can be unbounded. We refer to such a unbounded connected component as $\mathcal C_1$ and set
$$
\mathcal C_\lambda \,:=\,(\mathcal C_1^c\cap \Sigma_{\l_0^-})\cup R_\lambda(\mathcal C_1^c\cap \Sigma_{\l_0^-})
$$
If $u\equiv u_{\l_0^-}$ in $\mathcal C_1$ it is easy to see that, by symmetry, $\mathcal C_\lambda$ contains at least one connected component of $\mathbb{R}^N\setminus Z_u$. But this is not possible as it has been shown in \cite[Theorem 1.4]{DS1} and \cite[Lemma 5]{MMPS}.
 If else $u\equiv u_{\l_0^-}$ in $\mathcal C_2$ for some bounded component $\mathcal C_2$, then in this case we set
 $$
\mathcal C_\lambda\,:=\,\mathcal C_2\cup R_\lambda(\mathcal C_2)\,,
$$
 and also in this case, by symmetry,  $\mathcal C_\lambda$ would contain at least one connected component of $\mathbb{R}^N\setminus Z_u$ thus providing a contradiction. Resuming we just proved that
\begin{equation}\nonumber
u< u_{\l_0^-}\quad\text{in}\quad
\Sigma_{\l_0^-}\setminus{Z_u}\,.
\end{equation}
 Now, recalling that $Z_u$ is bounded by Theorem \ref{stimagrad},  we fix $\overline{R} >0$ in such a way that
 $$
 Z_u\subset B_{\overline{R}}(0)\,,
 $$
 and, for $\tau>0$, we let $Z_u^\tau$ be an open set containing $Z_u$ such that $\mathcal L (Z_u^\tau)<\tau$ (that exists since $\mathcal L (Z_u)=0$).
 Then, for  $\delta,\varepsilon, \overline{R},\tau >0$, we denote by
$$B_{\overline{R},\varepsilon}:= B^c_{\overline{R}}(0) \cap \Sigma_{\l_0^-+\varepsilon}, \ \
S_\delta^\varepsilon:= \left((\Sigma_{\l_0^-+\varepsilon}\setminus\Sigma_{\l_0^--\delta})\cap B_{\overline{R}}(0)\right)\cup (Z_u^\tau\cap \Sigma_{\l_0^--\delta}) ,$$$$  K_\delta:= \overline{B_{\overline{R}}(0)\cap \Sigma_{\l_0^--\delta}}\cap (Z_u^\tau)^c,$$
where $\delta\le \overline{\delta}$ so that $K_\delta$ is nonempty.
We underline that this construction gives $$\Sigma_{\l_0^- + \varepsilon} = B_{\overline{R},\varepsilon} \cup S_\delta^\varepsilon \cup K_\delta.$$

We also remark that, since $K_\delta$ is compact, then by the uniform continuity of $u$ and $u_\l$, for $\overline{\varepsilon}>0$ small enough one has that $u<u_{{\l_0^-}+\varepsilon}$ in $K_\delta$ for every $\varepsilon\le \overline{\varepsilon}$. Moreover we underline the existence of $\tilde{R_0}$ such that $u<u_{\l_0^- + \varepsilon}$ in $B_{\tilde{R_0}}(0_{\l_0^- + \varepsilon})\subset\Sigma_{\l_0^- + \varepsilon}$ for every $\varepsilon \le \overline{\varepsilon}$ and with $\tilde{R_0}$ independent of $\varepsilon$ as done in Step 1. \\

From now on, for $R>\overline{R}$, we consider $\eta\in C^\infty_0(B_{2R}(0))$ a cut-off function with $0\le \eta\le 1$, $\eta\equiv 1$ on $B_R(0)$ and $|\nabla \eta|\le \frac{2}{R}$.
Then, letting $\alpha>\max\{2, p\}$, we consider the following test functions
$$\varphi_{1, \l_0^- +\varepsilon}=\eta^\alpha u^{1-p}\(u^p-u_{\l_0^-+\varepsilon}^p\)^+\chi_{\Sigma_{\l_0^- + \varepsilon} },\qquad \varphi_{2, \l_0^- +\varepsilon}=\eta^\alpha u_{\l_0^-+\varepsilon}^{1-p}\(u^p-u_{\l_0^-+\varepsilon}^p\)^+\chi_{\Sigma_{\l_0^- +\varepsilon} },$$
and, analogously to Step $1$, $\psi_{\l_0^- +\varepsilon}:=(u^p-u_{\l_0^- +\varepsilon}^p)^+$ and by $\varphi_{\l_0^- +\varepsilon}:=(u-u_{\l_0^- +\varepsilon})^+$.
\\Let us take $\varphi_{1, \l_0^- +\varepsilon}$ as a test function in \eqref{pb}, $\varphi_ {2, \l_0^- +\varepsilon}$ in \eqref{pbriflessa} and, reasoning as in Step $1$, one yields to
\begin{equation}\label{stimastep2}
\begin{aligned}
&c_1 \int_{\hat B_{2R} \cap \{u\ge u_{\l_0^- +\varepsilon}\}}\eta^\alpha u^2\(|\n u|+|\n u_{\l_0^- +\varepsilon}|\)^{p-2}|\n\log u -\n\log u_{\l_0^- +\varepsilon}|^2
\\
&\le \int_{\hat B_{2R}\cap B_{\overline{R},\varepsilon}}(u^{p^*-p}-u_{\l_0^- +\varepsilon}^{p^*-p})\eta^\alpha\psi_\l + \int_{\hat B_{2R}\cap S_\delta^\varepsilon}(u^{p^*-p}-u_{\l_0^- +\varepsilon}^{p^*-p})\eta^\alpha\psi_{\l_0^- +\varepsilon}+\frac{C}{R^{\frac \beta p}}+ \frac{C}{R^\beta}.
\end{aligned}
\end{equation}
Here we have used once again the fact that $\frac{u_{\l_0^- +\varepsilon}}{u}\ge \tilde c$ for every $0\le \varepsilon\le \bar\varepsilon$ as to deduce \eqref{stimaulu}.\\
In order to estimate the first term on the right hand side of \eqref{stimastep2} we argue exactly as to estimate $I_4$ in \eqref{stima0} (taking into account Remark \ref{remngu}) where here $\overline{R}$ plays the role of $\l$ in Step $1$. Hence we get
\begin{equation*}\label{stima1step2}
\begin{aligned}
&\int_{\hat B_{2R}\cap B_{\overline{R},\varepsilon}}(u^{p^*-p}-u_{\l_0^- +\varepsilon}^{p^*-p})\eta^\alpha\psi_\l \le \frac{C}{R^\beta}
\\
&+ \frac{C}{\overline{R}^{\beta^*}}\int_{\hat B_{2R}\cap B_{\overline{R},\varepsilon}\cap \{u\ge u_{\l_0^- +\varepsilon}\}}\eta^\alpha u^2\(|\n u|+|\n u_{\l_0^- +\varepsilon}|\)^{p-2}|\n\log u -\n\log u_{\l_0^- +\varepsilon}|^2.
\end{aligned}
\end{equation*}
For the second term on the right hand side of \eqref{stimastep2} we reason as in Step $1$, getting
\begin{equation}\label{stima2step21}
\begin{aligned}
&\int_{\hat B_{2R}\cap S_\delta^\varepsilon}(u^{p^*-p}-u_{\l_0^- +\varepsilon}^{p^*-p})\eta^\alpha\psi_\l \le C_u\int_{\hat B_{2R}\cap S_\delta^\varepsilon\cap \{u\ge u_{\l_0^- +\varepsilon}\}}(\log u -\log u_{\l_0^- +\varepsilon})^2,
\end{aligned}
\end{equation}
where
$$C_u:=\begin{cases}
	\displaystyle \sup_{S_{\overline{\delta}}^{\overline{\varepsilon}}} u^{p^*-2} \quad \text{if } p^*\ge 2,\\
	\displaystyle \inf_{S_{\overline{\delta}}^{\overline{\varepsilon}}} u^{p^*-2} \quad \text{if } p^*< 2.
\end{cases}$$
Now we need to divide the estimate by the value of p; indeed if $p>2$ we apply a suitable weighted Poincar\'e inequality to the right hand side of \eqref{stima2step21} which can be found in Theorem $3.2$ of \cite{DS1}. Hence in this case one has
\begin{equation*}\label{stima2step2}
\begin{aligned}
&\int_{\hat B_{2R}\cap S_\delta^\varepsilon}(u^{p^*-p}-u_{\l_0^- +\varepsilon}^{p^*-p})\eta^\alpha\psi_\l
\\
&\le C^2_p(S_\delta^\varepsilon)C_u\int_{\hat B_{2R}\cap S_\delta^\varepsilon\cap \{u\ge u_{\l_0^- +\varepsilon}\}}|\n u|^{p-2}|\n\log u - \n \log u_{\l_0^- +\varepsilon}|^2
\\
&\le \frac{C^2_p(S_\delta^\varepsilon)C_u}{\displaystyle \inf_{S_{\overline{\delta}}^{\overline{\varepsilon}}}u^2}\int_{\hat B_{2R}\cap S_\delta^\varepsilon\cap \{u\ge u_{\l_0^- +\varepsilon}\}}u^2\(|\n u|+|\n u_{\l_0^- +\varepsilon}|\)^{p-2}|\n\log u - \n \log u_{\l_0^- +\varepsilon}|^2,
\end{aligned}
\end{equation*}
where $C_p(E)$ is the Poincar\'e constant which goes to zero as $|E|\to 0$. Otherwise if $1<p<2$ one can apply the classical Poincar\'e inequality in order to deduce
\begin{equation*}\label{stima2bisstep2}
\begin{aligned}
&\int_{\hat B_{2R}\cap S_\delta^\varepsilon}(u^{p^*-p}-u_{\l_0^- +\varepsilon}^{p^*-p})\eta^\alpha\psi_\l
\\
&\le C^2_p(S_\delta^\varepsilon)C_u\int_{\hat B_{2R}\cap S_\delta^\varepsilon\cap \{u\ge u_{\l_0^- +\varepsilon}\}}|\n\log u - \n \log u_{\l_0^- +\varepsilon}|^2
\\
&\le \frac{CC^2_p(S_\delta^\varepsilon)C_u}{\displaystyle \inf_{S_{\overline{\delta}}^{\overline{\varepsilon}}}u^2}\int_{\hat B_{2R}\cap S_\delta^\varepsilon\cap \{u\ge u_{\l_0^- +\varepsilon}\}}u^2\(|\n u|+|\n u_{\l_0^- +\varepsilon}|\)^{p-2}|\n\log u - \n \log u_{\l_0^- +\varepsilon}|^2,
\end{aligned}
\end{equation*}
which can be deduced since in $\Sigma_{\l_0^- + \varepsilon} \setminus B_{\tilde{R}_0}(0_{\l_0^- + \varepsilon})$ one has that
	$$\(|\n u|+|\n u_\l|\)^{2-p}\le C,$$
for some constant $C$ which does not depend on $\varepsilon\le \overline{\varepsilon}$.
Hence in both cases one has that
\begin{equation*}\label{stima3step2}
\begin{aligned}
&c_1\displaystyle \int_{\hat B_{2R}\cap \{u\ge u_{\l_0^- +\varepsilon}\}}\eta^\alpha u^2\(|\n u|+|\n u_{\l_0^- +\varepsilon}|\)^{p-2}|\n\log u -\n\log u_{\l_0^- +\varepsilon}|^2
\\
&  \le  \frac{C}{R^{\frac \beta p}}+\frac{C}{R^\beta}+\frac{C}{\bar R^{\beta^*}}\displaystyle \int_{(\hat B_{2R}\cap B_{\bar R, \varepsilon})\cap \{u\ge u_{\l_0^- +\varepsilon}\}}\eta^\alpha u^2\(|\n u|+|\n u_{\l_0^- +\varepsilon}|\)^{p-2}|\n\log u -\n\log u_{\l_0^- +\varepsilon}|^2\\ &+\frac{CC^2_p(S_\delta^\varepsilon)C_u}{\displaystyle \inf_{S_{\overline{\delta}}^{\overline{\varepsilon}}}u^2}\int_{S_\delta^\varepsilon\cap \{u\ge u_{\l_0^- +\varepsilon}\}}\(|\n u|+|\n u_{\l_0^- +\varepsilon}|\)^{p-2}|\n\log u - \n \log u_{\l_0^- +\varepsilon}|^2.
\end{aligned}
\end{equation*}

Now we take care of the variable parameters $\bar R, \delta, \bar\varepsilon$. First we fix $\bar R$ large such that  $$\frac{C}{c_1 \bar R^{\beta^*}}<1.$$ Then, since $C_p^2(\Omega)$ goes to zero if the Lebesgue measure of $\Omega$ goes to zero, we choose $\delta,\bar\varepsilon,\tau$ small so that $$\frac{CC^2_p(S_\delta^\varepsilon)C_u}{c_1\displaystyle \inf_{S_{\overline{\delta}}^{\overline{\varepsilon}}}u^2}<1$$ for every $0\le \varepsilon \le \bar\varepsilon$.
Hence it follows that
\begin{equation*}\label{stima4step2}
\begin{aligned}
&\int_{\hat B_{2R}\cap \{u\ge u_{\l_0^- +\varepsilon}\}}u^2\(|\n u|+|\n u_{\l_0^- +\varepsilon}|\)^{p-2}|\n\log u -\n\log u_{\l_0^- +\varepsilon}|^2
\le \frac{C}{R^{\frac \beta p}}+\frac{C}{R^\beta}
\end{aligned}
\end{equation*}

getting again (as $R\to+\infty$)
\begin{equation*}\label{stima6step2}
\begin{aligned}
&\int_{\Sigma_{\l_0^- + \varepsilon} \cap \{u\ge u_{\l_0^- +\varepsilon}\}}u^2\(|\n u|+|\n u_{\l_0^- +\varepsilon}|\)^{p-2}|\n\log u -\n\log u_{\l_0^- +\varepsilon}|^2 = 0,
\end{aligned}
\end{equation*}
which gives that $u\le u_{\l_0^- +\varepsilon}$ in $\Sigma_{\l_0^- + \varepsilon}$ which contradicts the definition of $\l_0^-$. This proves that $\l_0^-=0$. In an analogous way we deduce that $\l_0^+=0$, which gives the symmetry of $u$ along the $e_1$-direction. Repeating the same arguments in the remaining $N-1$ linearly independent directions of $\R$ then one deduces that $u$ is symmetric about the origin and that is a radially decreasing function.
\end{proof}

\end{document}